\documentclass{article}
\usepackage[utf8]{inputenc} 
\usepackage[T1]{fontenc}    
\usepackage{hyperref}       
\usepackage{url}            
\newcommand{\tarc}{\mbox{\large$\frown$}}
\newcommand{\arc}[1]{\stackrel{\tarc}{#1}}
\usepackage{booktabs}       
\usepackage{amsfonts}       
\usepackage{nicefrac}       
\usepackage{graphicx,subfigure}
\graphicspath{ {images/} }
\usepackage{microtype}      
\usepackage{lipsum}		
\usepackage{float}
\usepackage{graphicx}
\usepackage{natbib}
\usepackage{doi}
\usepackage{amssymb}
\usepackage{amsmath}
\DeclareUnicodeCharacter{2212}{-}
\setcitestyle{numbers,square}
\let\oldcdot\cdot
\usepackage{breqn}
\let\cdot\oldcdot

\title{ Solution to the Riemann Hypothesis from geometric analysis of component series functions in the functional equation of zeta}

\author{Jeet Kumar Gaur\\\small PhD in Mechanical Eng. Dept. \hfill \\\small Indian Institute of Science, Bangalore, India \\\small email id: jeetgaur@iisc.ac.in}

\hypersetup{
pdftitle={A B C D},
pdfsubject={Riemann hypothesis},
pdfauthor={JKG},
pdfkeywords={Riemann Hypothesis, Proof},
}
\begin{document}
\maketitle
\begin{abstract}
  This paper presents a new approach towards the Riemann Hypothesis. On iterative expansion of integration term in functional equation of the Riemann zeta function we get sum of two series function. At the `non- trivial' zeros of zeta function, value of the series is zero. Thus, Riemann hypothesis is false if that happens for an `s' off the line $\Re(s)=1/2$ ( the critical line). This series has two components $f(s)$ and $f(1-s)$. For the hypothesis to be false one component is additive inverse of the other. From geometric analysis of spiral geometry representing the component series functions $f(s)$ and $f(1-s)$ on complex plane we find by contradiction that they cannot be each other's additive inverse for any $s$, off the critical line. Thus, proving truth of the hypothesis.   
\end{abstract}
\section{Introduction}
The Riemann Hypothesis has become one of the most centralized problems in mathematics today. It has to do with the position of `non-trivial' zeros of the zeta function that encrypt information about prime numbers. The hypothesis has its roots predominantly in number theory. However, over time the ways to approach the problem have diversified with new evidence suggesting connections with seemingly unrelated fields like random matrices, chaos theory and quantum physics \cite{1}\cite{QM}. In over  160 years of period, from the time it was first formulated in 1859, Riemann hypothesis has been attempted by several great mathematical minds. But, yet it has remained unsolved.   Solving Riemann hypothesis would also verify hundreds of theorems and solutions which are based on the assumption of it's truth.\\\\
Owing to the works of Turing, Hardy, Siegel, Chebyshev and others, better algorithms to numerically evaluate the zeta function have been developed \cite{comp}. As of now the first 10 billion zeros have been numerically verified to follow the hypothesis without a single exception \cite{comp}. Although computational pieces of evidence speak in favour of the hypothesis, an analytic approach to the hypothesis is necessary to prove it for each and every zero.\\\\
The Riemann zeta function, $\zeta(s)$ is defined by a Dirichlet series as \\
\begin{equation*}
\zeta(s) = \sum\limits_{n=1}^{\infty} \frac{1}{n^{s}}= 1 + \frac{1}{2^{s}} + \frac{1}{3^s} + \frac{1}{4^s} + ...
\end{equation*}\\
 which is only convergent in the  region  $\Re(s)>1$, \cite{Zseries}.
\\\\
The zeta function was first known to the Swiss Mathematician Leonhard Euler. In 1737, Euler gave the famous form of the zeta function as a product of functions of prime numbers \cite{PO99}\cite{P100}.\\
\begin{equation*}
\zeta(s) = \prod\limits_{p} \left(1-\frac{1}{p^{s}}\right)^{-1}.
\end{equation*}
The above equation also represents the secret association of prime numbers with  the zeta function. However, it was Bernhard Riemann who first studied the complex form of $\zeta(s)$ and also found a functional equation that analytically continues it as a meromorphic function. He found that the concrete prediction for the occurrence of prime numbers depends on the placement of `non-trivial'( complex) zeros of this function. If the Riemann hypothesis is true then all these zeros lie collinearly. The values of the function on negative integers is associated to the corresponding Bernoulli number. The Bernoulli numbers can be obtained from Euler-Maclauren series(1735) \cite{Bern}\cite{Wone}. The $n$th Bernoulli number $B^{+} _{n}$ can be defined as\\
\begin{equation*}
\begin{aligned}
{B_{n}^{+} = − n*\zeta(1-n)},
 \end{aligned}
\end{equation*}\\
for all $n \geq 1$. Every odd Bernoulli number $B^{+} _{2n+1}$ is zero so, the value of zeta function at every negative even integer. These are the `trivial' (real) zeros of zeta function\\\\
In 1850's Bernhard Riemann, a theology student in Germany turned back to his first love mathematics and became a student of Friedrich Gauss. In his Legendary  1859 paper titled "Ueber die Anzahl der Primzahlem unter eine gegebener Gr${\ddot{o}}$sse" that is "On the Number of Primes Less Than a Given Magnitude", Riemann conjectured his hypothesis. He presented the unnatural connection between the position of zeros of the zeta function with order of the error term in prime counting functions such as conjectured by Gauss. The Prime Number Theorem( PNT) states that the number of primes less than or equal to $n$ approaches $n/\log{n}$ as $n$ approaches infinity \cite{PNT}. \\
\begin{equation*}
\lim\limits_{n\to\infty}\frac{\pi(n)}{n/\log(n)}=\;1.
\end{equation*}
Gauss shared his prime counting formula with Dirichlet who came up with a more accurate counting function, the $Li$ function.\\
\begin{equation*}
Li(x):=\int\limits_{2}^{x}\frac{dx}{\log(x)}.
\end{equation*}\\
Although $Li(x)$ and $\pi(x)$ asymptotically become more and more accurate as $n$ increases but none of them is a true prime counting function. The Riemann hypothesis is connected to the order of the error term. If the Riemann hypothesis is true, there is an error of order $\sqrt{x}$ in counts of primes as given by the $Li$ function \cite{error}. That is for arbitrary $\epsilon>0$, 
\begin{equation*}
\pi_{true}(x)=\int\limits_{2}^{x}\frac{dx}{\log(x)} + O(x^{\frac{1}{2}+\epsilon}).
\end{equation*}\\
There are several other equivalent statements for RH made using it's connection with functions from number theory. For example, the $Liouville function$, $\lambda(n)$ which is defined by\\
\begin{equation*}
\lambda(n)\;=\;(-1)^{\omega(n)},\\
\end{equation*}
where $\omega(n)$ is the number of prime factors, counted with multiplicity. Thus, $\omega(6)= 2$ and $\lambda(6)=\;1$,\\
$\omega(9)= 2$ and $\lambda(9)=1$ and\\
$\omega(p)=1$ and $\lambda(p)=-1$ for all primes $p$.\\
In his doctoral thesis, Laudau (1899) showed that RH is equivalent to the statement\\
\begin{equation*}
\lim\limits_{n\to\infty} \frac{\lambda(1)\;+\;\lambda(2)\;+\;\lambda(3)\;+\;.\;.\;.+\;\lambda(n)}{n^{\frac{1}{2}+\epsilon}}\;=\;0,
\end{equation*}
for all $\epsilon > 0$  \cite{56}.
The theory of RH is vast and connections to seemingly unrelated fields have been found. One such example is the connection to Random Matrix Theory as explored by Montgomery. In 1972, Hugh Montgomery investigated the spacing between the zeros of the zeta function. This was just after the spectral interpretation of RH provided by Hilbert-P\'olya conjecture \cite{1}. Assuming RH to be true, Montgomery conjectured that the correlation for the zeros of the zeta function is\\
\begin{equation*}
1-\frac{\sin^{2}(\pi x)}{(\pi x)
^2}.
\end{equation*}
This is the pair correlation Conjecture given by Montgomery. Later Freeman Dyson identified that the above correlation function is the same as that for pairs of eigenvalues of random Hermitian matrices with Gaussian measure used by physicists. This gave birth to the Gaussian Unitary Ensemble( GUE) conjecture\cite{JBC}.\\
Till now several equivalent statements and extensions of RH have been formulated. New findings on association of zeta function with other concepts keep redefining the hypothesis and encourage newer approach to the problem. Despite of having acquired several pieces of information about the zeta function, the problem has remained unsolved till date. However, famous attempts to solve RH are worth mentioning.\\\\
In 1885,  T. J Stieltjes claimed to have solved the Mertens' conjecture which implied the Riemann hypothesis. The results were proven wrong after his death as Mertens’ conjecture was proven false by Andrew Odlyzko and Herman J.J. te Riele in 1985 \cite{stlj}.\\\\
In the meanwhile, German-born American mathematician Hans Rademacher came up with a proof of the falsehood of Riemann Hypothesis. He reported his proof in \textit{Times magazine} in 1945 but it was found erroneous. Legend has it that Rademacher was terribly embarrassed by the ordeal and it was well known that no one was to mention the words, \textit{Riemann Hypothesis} in his presence\cite{Rade}.\\\\ In 1959 the famous mathematician John Nash presented a proof for the Riemann hypothesis in a lecture at Columbia University. The attendees had high expectations but were disappointed as the lecture didn't make any sense. This was later found to be due to his struggle with schizophrenia \cite{nash}.\\\\
In 2004, Louis de Bourcia claimed to have the proof of Riemann hypothesis which remained under skepticism for a long time. Later a counterexample from a 1998 paper by one of his student mathematicians Xian-Jin Li and Brian Conrey showed his proof was false. Li, in July 2008 reported the proof of RH but it was soon retracted after  mathematicians pointed out a crucial flaw\cite{Lin}.\\\\
The most recent famous claim came in September, 2018 from British-Lebenese mathematician Sir Micheal Atiyah. He used Todd functions and mentioned the Fine Structure constant a fundamental physical constant in his proof by contradiction \cite{atiy}. The Fine Structure constant also called the Sommerfeld's constant is a dimensionless constant that equals $1/137$. A concept entirely from physics. His proof was based on the works of leading 20th century mathematicians, John von Neumann and Friedrich Hirzebruch. His proof remained under skepticism and to be rigorously verified. Although few mathematicians have declared it flawed \cite{new}. \\  
\section{Expansion of functional equation of zeta function into a series} 
 The Dirichlet series form of the zeta function is only valid for $\Re(s)>1$. Riemann gave the functional equation of Zeta function that continued it as a meromorphic function with a simple pole at 1 \cite{1}.

\begin{equation*}
\zeta(s) :=  \frac{\pi^{s/2}}{\Gamma(s/2)}\Bigg\{\frac{1}{s(s-1)} + \int\limits_{1}^{\infty}(x^{s/2-1} + x^{-s/2-1/2}).\left(\frac{\vartheta(x)-1}{2}\right)dx\Bigg\}
\end{equation*}

defined for Re(s) > 0, where
\begin{equation*}
\vartheta(x) := \sum\limits_{n=-\infty}^{\infty}e^{-n^2{\pi}x} = 1 + 2\sum\limits_{n=1}^{\infty}e^{-n^2{\pi}x} 
\end{equation*}
is the Jacobi-theta function.\\\\
From analytic continuation of zeta function, one also gets the  Riemann reflection formula 
 \begin{equation*}
   \xi(s)=\xi(1-s).
 \end{equation*}
Where $\xi(s)$ is an entire function and is defined as\\
\begin{equation*}
\xi(s)=\frac{s}{2}(s-1)\pi^{-\frac{s}{2}}\Gamma(\frac{s}{2})\zeta(s).    
\end{equation*}
 In his extensive work on zeta function in 1859, Riemann also gave his famous hypothesis about the position of non-trivial zeroes of the zeta function. The Riemann hypothesis is essentially that all the values of $s$ in the critical strip for which $\zeta(s) = 0$, lie on the critical line ($\Re(s) =\frac{1}{2}$) \cite{1}.\\
Let us define:
\begin{equation*}
    {\textstyle\Omega(s)=\frac{1}{s(s-1)} +  \int\limits_{1}^{\infty}(x^{\frac{s}{2} -1}+x^{-\frac{s}{2}-\frac{1}{2}}).\left( \frac{\vartheta\left( x\right)  -1}{2} \right) dx} 
\end{equation*}\\
defined in $(0 < \Re(s) < 1)$.\\
On expanding the integral term in $\Omega(s)$, by ‘integration by parts’ taking $x^{\frac{s}{2}-1}+x^{-\frac{s}{2}-\frac{1}{2}}$ as the first function and the Jacobi-theta function term as the
second,
\begin{equation*}
\begin{aligned}
\textstyle\Omega(s)&=\frac{1}{s(s-1)} + \sum\limits_{n=1}^{\infty} \left(0 -\frac{2}{s}-\frac{2}{1-s}\right)e^{-n^2\pi}\\
&-\int\limits_{1}^{\infty}\left(  \frac{2 x^{\frac{s}{2}}}{s} + \frac{2 x^{-\frac{s}{2}+ \frac{1}{2}}}{1-s}\right) \sum\limits_{n=1}^{\infty} e^{-n^2\pi x}(-n^{2}\pi) dx
\end{aligned}
\end{equation*}
Repeating integration ‘by parts’ again and again, taking the Jacobi-theta function part as the second function, one gets a series function:

\begin{equation*}
\begin{aligned}
\Omega(s) = &\frac{1}{s(s-1)} + \sum\limits_{n=1}^{\infty} \left(-\frac{2}{s}-\frac{2}{1-s}\right)e^{-n^2\pi}\\
&+ \sum\limits_{n=1}^{\infty} \left(-\frac{1}{\frac{s}{2}(\frac{s}{2} +1)} -\frac{1}{(\frac{1-s}{2})(\frac{1-s}{2}+1)}\right)e^{-n^2\pi}.n^2\pi \\
&+ \sum\limits_{n=1}^{\infty} \left(-\frac{1}{\frac{s}{2}(\frac{s}{2} +1)}-\frac{1}{(\frac{1-s}{2})(\frac{1-s}{2}+1)}\right)e^{-n^2\pi}.(n^2\pi)^2 + ...
\end{aligned}
\end{equation*}

We define:
$$ f(s)= -\frac{1}{s} +\sum\limits_{n=1}^{\infty} e^{-n^2\pi} \left(
-\frac{1}{\left(\frac{s}{2}\right)^{\overline {1}}} - \frac{n^2\pi}{\left(\frac{s}{2}\right)^{\overline{2}}}- \frac{(n^2\pi)^2}{\left(\frac{s}{2}\right)^{\overline{3}}}- \frac{(n^2\pi)^3}{\left(\frac{s}{2}\right)^{\overline {4}}}- ...\right)$$

where $(\frac{s}{2})^{\overline{n}}=\frac{s}{2}(\frac{s}{2}+1)(\frac{s}{2}+2)...(\frac{s}{2}+n-1)$ is the rising factorial power function.\\\\
So, 
\begin{equation*}
\Omega(s) = f(s) + f(1-s),
\end{equation*}
and
\begin{equation*}
\textstyle f(s)=-\left(\frac{1}{s}+J(s)\right).
\end{equation*}\\
Where, \begin{equation*}
J(s):= \sum\limits_{m=1}^{\infty}\frac{\sum\limits_{n=1}^{\infty} e^{-n^2\pi}\left(n^2\pi\right)^{m-1}}{\left(\frac{s}{2}\right)^{\overline{m}}}.
\end{equation*}\\
\section{Proof of the truth of RH}
The zeros of $\Omega(s)$ are the non-trivial zeros of the zeta function. 
$\Omega(s)$ can be equal to zero in either of the two cases:
\begin{itemize}
    \item if $f(s)$ and $f(1-s)$ are individually zero.
    \\\\Or
    \item if $f(s)=-f(1-s)$ and none equals zero.
    \end{itemize}
    Therefore, proving that $f(s)$ is not the additive inverse of $f(1-s)$ anywhere off the critical line is sufficient to proof the Riemann Hypothesis.\\\\
    In the next section deductions about the allowed values of the series function $J(s)$ are made by geometric analysis of simpler series functions and $J(s)$. From the allowed values of $J(s)$, the existence of the first condition for $\Omega(s)=0$ is obtained.      
\subsection{ \texorpdfstring{Geometry of series $J(s)$}{lg} on the complex plane}
Plotting $J(s)$ term by term on the complex plane for a given \textit{s}, leads us to value of the series. In the present section we will find that $J(s)$ converges as an anticlockwise inward spiral starting from origin, for values of \textit{s} in the lower half region of the critical strip.\\ 
Let's define $\alpha \in (0,\frac{1}{2})$.
Thus, $s=\frac{1}{2} \pm \alpha - i \vert t \vert$ represents any point in the lower region of the critical strip to the real axis (of course, excluding the critical line $\Re(s)=\frac{1}{2}$). This is our region of interest because if Riemann hypothesis is false then $\Omega(s)$ is zero in this region.\\
Now, for $s=\frac{1}{2} \pm \alpha - i \vert t \vert$, the argument of the $m^{th}$ term of $J(s)$ is given by:
\begin{equation*}
\begin{aligned}
\theta_{m,J} & = arg\left(\frac{1}{\frac{s}{2}}\right)^{\overline{m}}= -arg \left(\frac{s}{2}\right)^{\overline{m}}  \\ 
&=\arctan\left(\frac{\frac{\vert t \vert}{2}}{\frac{\Re(s)}{2}}\right)+\arctan\left(\frac{\frac{\vert t \vert}{2}}{\frac{\Re(s)}{2}+1}\right)+...+\arctan\left(\frac{\frac{\vert t \vert}{2}}{\frac{\Re(s)}{2}+m-1}\right)
\end{aligned}
\end{equation*}
So, $\theta_{m,J} < m \theta$, where $\theta$ is $\theta_{1,J}$ (argument of the first term OA).
$0<\theta_{m+1,J}-\theta_{m,J} = \arctan\left(\frac{\frac{\vert t \vert}{2}}{\frac{\Re(s)}{2}+m}\right)< \frac{\pi}{2}$ i.e., is positive and acute $\forall\;\;m$  .
Thus, as $m$ increases, the difference between the arguments of consecutive terms of $J(s)$ decreases.

Since, $J(s)$ is a component function in the functional equation of zeta function, it is also analytic throughout the critical strip. A proof of its convergence is provided in the appendix. Form its convergence we have,  
$$ \frac{\vert T_{m+1,J}\vert}{\vert T_{m,J}\vert}= \frac{\sum\limits_{n=1}^{\infty} e^{-n^2\pi}\left(n^2\pi\right)^{m}}{\vert\frac{s}{2} + m\vert \sum\limits_{n=1}^{\infty} e^{-n^2\pi}\left(n^2\pi\right)^{m-1}}<1 \;\;\; \forall\;\; m
$$
\begin{figure}[!htb]
\subfigure[$\frac{\vert T_{m+1}\vert}{\vert T_{m}\vert} = 1,\; \theta_{m}=m\theta\; and \;\frac{2\pi}{\theta} \in\; \mathbb{N}^{+}$]
{\includegraphics[scale=0.77]{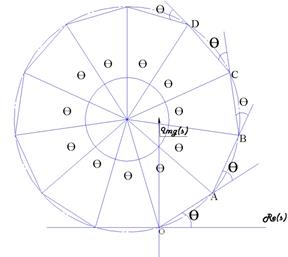}}
\subfigure[$\frac{\vert T_{m+1}\vert}{\vert T_{m}\vert} = 1,\;  \theta_{m}=m\theta\; and \;\frac{2\pi}{\theta} \notin\; \mathbb{N}^{+}$]
{\includegraphics[scale=0.78]{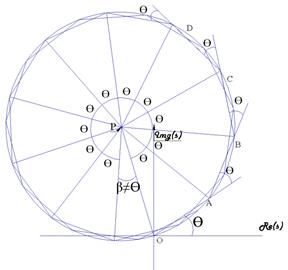}}
\caption{Different cases of series $U(s)$}
\label{fig:image1}
\end{figure}
Now, let's try to plot $J(s)$, by plotting some simpler series.\\
Take a series,\\
$$ U(s)=\sum\limits_{m=1}^{\infty} T_{m}(s) = \sum\limits_{m=1}^{\infty}{\vert T_{m}(s) \vert} e^{i\theta_{m}}$$\\
such that $$\vert T_{m+1}(s) \vert=\vert T_{m}(s) \vert=T$$

Therefore, we have
$$  OA = T_{1}(s)=\vert T_{1}(s)\vert e^{i\theta}=T e^{i\theta}$$
$$AB = \vert T_2(s) \vert e^{2i\theta}= T e^{2i\theta}$$
$$BC = T_3(s)=\vert T_{3}(s) \vert e^{3i\theta}= T e^{3i\theta}$$ 
\\
Fig. \ref{fig:image1} shows U(s) series with $ \frac{\vert T_{m+1}\vert}{\vert T_{m}\vert}=1$, $ \theta_{m}=m\theta $, i.e. $\frac{T_{m+1}(s)}{T_{m}(s)}= e^{i\theta}$. Here, either of the following two cases will occur:

\begin{itemize}
  \item A uniform polygon of sides $\frac{2\pi}{\theta}$ : the whole series plots as a single polygon of $\frac{2\pi}{\theta}$ sides (fig. a).
  \item if $\theta$ does not divide $2\pi$ i.e., the circle into $n$ equal sectors and we have a geometry as shown (fig. b) with consecutive cycles going out of phase and there exists $\beta: 0 < \beta < \theta$, and $\beta + (n-1)\theta = 2\pi$, which is the angle formed by the last sector at the center of the polygon on completion of $2\pi$.\\
  \begin{figure}[!hb]
\centering     
\includegraphics[scale=.3]{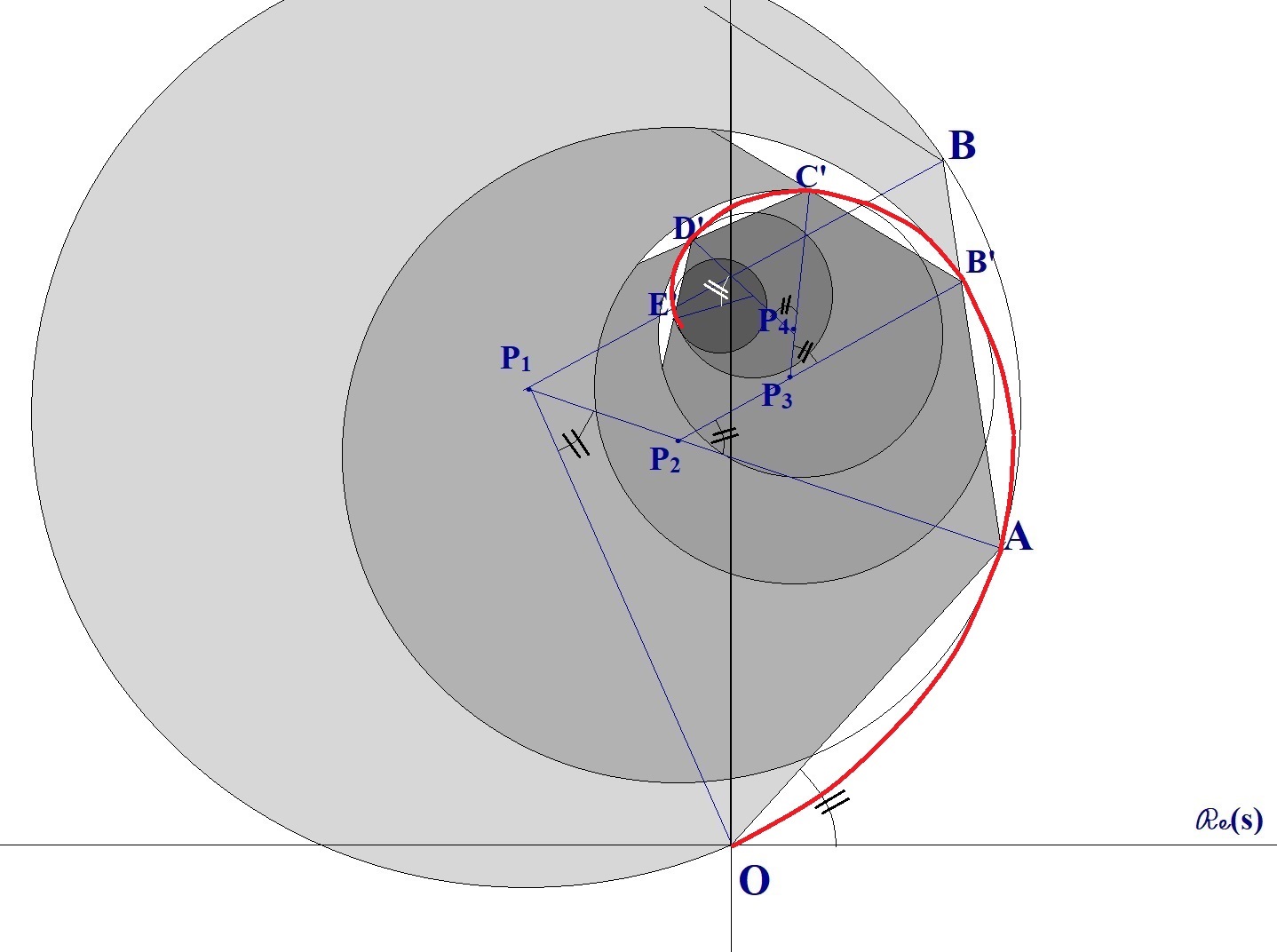}
\caption{$\mu(s): \frac{\vert T_{m+1}\vert}{\vert T_{m}\vert} < 1$, $ \theta_{m}= m \theta$}
\label{fig:image5}
\end{figure}
  \end{itemize}
  But, in none of the cases above, does the series converge.\
Let's define a new series function $\mu(s)=\sum\limits_{m=1}^{\infty} T_{m,\mu}(s)$, such that $\frac{\vert T_{m+1,\mu}\vert}{\vert T_{m,\mu}\vert}<1$ and, $\theta_{m,\mu}=\theta_{m}=m\theta$, for $\mu(s)$,
 where $\theta$ is acute and positive.\\
 
As shown in fig. \ref{fig:image5}, we get an inwardly spiraling geometry going anticlockwise from the origin. The value of  $\mu(s)$ always lie inside the circle passing through points O, A, B, C, D, .... 
Fig. \ref{fig:image5} shows that any series with the attributes of $\mu(s)$ necessarily lies in the circle with centre $P_{1}$ and radius $AP_1=BP_1$. Take the first term $OA$ as in the $U(s)$ series.\\
Since, $$\vert T_{2,\mu}(s)\vert <\vert T_{1,\mu}(s)\vert$$ thus B' lies between points A and B on chord $AB$.
Thus, the point B' lies inside circle with centre $P_1$. Now, draw a circle with AB' as its chord such that it projects the same angle at it's center $(P_2)$ as chord AB projects at the center ($P_1$) of outer circle. This new circle with centre $P_2$ and radius $B'P_2=AP_2\;< AP_1$ is drawn by moving center $P_{1}$ along $AP_1$  to $P_2$ such that it's tangential at point A to the outer circle and $B'P_2$ is parallel to $BP_1$. Using the similarity of $\bigtriangleup AP_1B$ with $\bigtriangleup AP_2B'$, the length $AP_2$ follows $\frac{AP_2}{AP_1} =\frac{AB'}{AB}\;<\;1$. Since, circle with centre $P_2$ touches the outer circle only tangentially at point A, it's area always lie inside that of the outer circle. Same treatment can be done for the third term  $T_{3,\mu}(s)$ represented by $B'C'$ which is less than $BC$ by generating new circle tangential to the previous one at point B'. Thus circle with center $P_3$ and radius $C'P_3=B'P_3$\;($<\;B'P_2$ since  $\vert T_{3,\mu}(s)\vert<\vert T_{2,\mu}(s)\vert$) is generated by moving poimt $P_2$ along $B'P_2$ to point $P_3$. Again, area of the circle with center $P_3$ and radius $C'P_3$  lies inside the circle with center $P_2$. If $C_n$ represents the area under the $nth$ circle with center $P_{n}$ . We have,
$$C_n\subset\;C_{n-1}\subset\;...C_{2}\subset\;C_{1} $$. The value of the series $\mu(s)$ lies in the $nth$ circle as $n$ approaches infinity. The area of  circles with center $P_2, P_3, ..,P_n$ lie inside area $C_1$. This shows that the value of all series with attributes of $\mu(s)$ must lie inside the circle with center $P_1$ circumscribing the polygon formed for the case $T_{m+1}(s) = T_m(s)$ as represented by $U(s)$ above. The red colored spiral is formed by the arcs
$\arc{OA}$, $\arc{AB'}$, $\arc{B'C'}$, $\arc{C'D'}$ ... and converges at the same point as the series does.  \\
\textbf{The rate of convergence depends on the ratio of the consecutive terms. The  lower the ratio of consecutive terms $\frac{\vert T_{m+1}\vert}{\vert T_{m}\vert}$, the higher the rate of convergence of the series $\mu(s)$}.
Thus, on choosing the first term in $\mu(s)$ same as for $U(s)$, it converges at a
point inside the circle circumscribing the geometry formed for the
$\frac{\vert T_{m+1}\vert}{\vert T_{m}\vert}=1$
cases above. The grey to black circular regions represent necessary region in which the value converges after first few term.
\\\\
Now, coming back to our original series\\ $$J(s)=\sum\limits_{m=1}^{\infty} T_{m}(s)=\sum\limits_{m=1}^{\infty}\vert T_{m}(s)\vert e^{i \theta_{m,J}}$$.\\
where,  $ \theta_{m,J} = arg\left(\frac{1}{\frac{s}{2}}\right)^{\overline{m}}= -arg \left(\frac{s}{2}\right)^{\overline{m}}$ 
\\ 
\begin{figure}[!htb]
\centering     
\includegraphics[scale=.6]{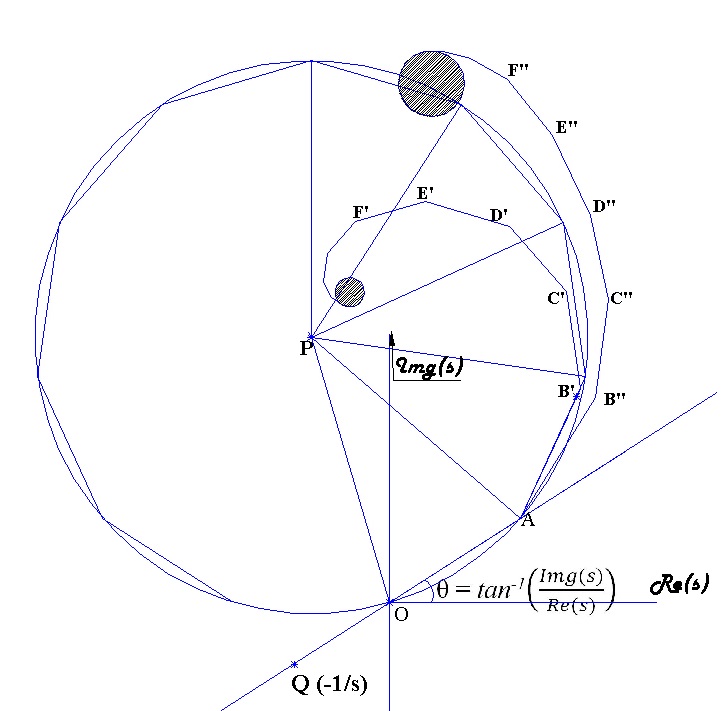}
\caption{$J(s): \frac{\vert T_{m+1}\vert}{\vert T_{m}\vert} < 1$, $ \theta_{m,J}< m \theta$ and $ \mu(s)$}
\label{fig:image4}
\end{figure}

We exploit the fact that for both $\mu(s)$ and $J(s)$ we have, $\frac{\vert T_{m+1}\vert}{\vert T_{m}\vert}<1$ and
draw $\mu(s)$ with terms having magnitude as that for $J(s)$ but $\theta_{m}=m \theta$. This corresponds to the geometry OAB'C'D'E'... in fig. $(3)$ below. 
 Now, we just have to rotate each term (starting from the second term) relative to the preceding one, clockwise. This rotation or unfastening is fixed for each term.  
So, from the last case we fix OA and rotate AB’ by
($\theta_{2,J} - \theta$) clockwise about A, to new position AB”, so that the acute angle
between AB” and line colinear with OA is now $\theta_{2,J}$.
Consecutively, rotate C’B”
clockwise about B” by
$\theta_{3,J} - \theta$
to new position B”C” and D’C” about C” by
$\theta_{4,J} - \theta$
and so on.
As shown in fig. \ref{fig:image4}, $J(s)$ is represented by OAB”C”D”E”... which is unfastened
form of the spiral OAB’C’D’E’.... Thus the rate of convergence is reduced. The unfastening is not done uniformly to all line segments representing series terms.
Rather, the unfastening monotonously decreases over terms, from
$arctan\left(\frac{\vert t \vert/2}{\Re(s)/2+1}\right)$
 at m= 2 (since
OA is fixed) to 0 as $ m\rightarrow \infty $.\\

In fig. \ref{fig:image4}, the spiral formed in fig. \ref{fig:image5} is unfastened so that, argument of $mth$ term is 
$\theta_{m,J} (<  m\theta) $
and depending on the difference 
$\theta_{m,J}-m\theta,$ we get the extent of unfastening( clockwise rotation) of each term of the spiral. Nonetheless this unfastened spiral also has a monotonous change in  curvature from its centre to it's farthest point at the origin, $O$. This is true for all series with $ \vert T_{m+1}\vert < \vert T_{m} \vert$ and $\theta_{m+1}>\theta_{m}$. The same applies to $J(s)$. 
 
Using the above deductions about $f(s)$ and $J(s)$, we look at the conditions for $\Omega(s)$ to be zero for any $\alpha\in\left(0,\frac{1}{2}\right)$.  

\subsection{Conditions for the Falsification of the Riemann Hypothesis}
\subsubsection{Case-1: RH is false if\texorpdfstring{ $f(s)=f(1-s)=0$, for $s=\frac{1}{2} \pm \alpha -i\vert t\vert$}{Lg}}
The condition $f(s)=0$ implies $J(s)=-\frac{1}{s}$. For any monotonously converging spiral, the distance from the centre of the spiral to a point on the spiral changes monotonously along the spiral. 

 Let point Q represent the value $-\frac{1}{s}$ on the complex plane. 
$$OA= \vert\; T_{1,J}(s)\;\vert =\;\frac{2\sum\limits_{n=1}^{\infty} e^{-n^2\pi}}{\vert s\vert}= \frac{\frac{\pi^{\frac{1}{4}}}{\Gamma(\frac{3}{4})}-1}{\vert s\vert}=\frac{0.0864348}{\vert s\vert}$$
\begin{figure}[!htb]
\centering     
\includegraphics[scale=.6]{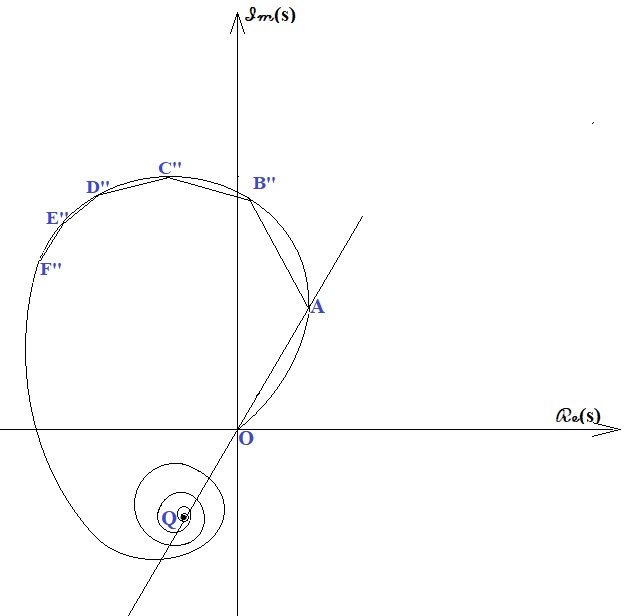}
\caption{Value of $J(s)$ lies col-linearly with chord OA if J(s) =-$\frac{1}{s}$}
\label{fig:image6}
\end{figure}
If  (see fig. \ref{fig:image6}) $J(s)=-\frac{1}{s}$ then, the the centre of spiral $OAB"C"D"E"...$ coincides point Q. Since, $J(s)$ is a monotonously converging spiral starting anticlockwise from first term $OA$ as discussed in previous section. Thus, $$...< QD"< QC"< QB"< QA < QO$$ as distance from Q to points on the spiral change monotonously as you move along the spiral. But, one finds that distance $QO=\frac{1}{\vert s\vert}$ and $QA= OA+OQ =\frac{1.0864348}{\vert s\vert}$. Thus, $$QO < QA.$$ Thus by contradiction
$\frac{1}{s} + J (s) \neq 0 $,
$$\Rightarrow f(s) \neq 0.$$

Also, since $(f(s))^{*} = f(s^{*})$ thus, $f(s)\neq O $   for any $s$, off the critical line, $\Re(s)=\frac{1}{2}$. 

Thus, no such case exists to falsify RH.

\subsubsection{Case-2: RH is false if \texorpdfstring{$f(s)=-f(1-s)$}{Lg}}
 For this case we first show that the magnitude of our complex series function, $|f(s)|$ always decreases as the real part of $s$ changes from zero to one keeping the imaginary part of $s$ constant. Since, showing that for $\Re(s)$ less than $\frac{1}{2}$, the magnitude $|f(1-s*)|$ is always less than $|f(s)|$, is sufficient to show that $f(s)\neq - f(1-s)$ anywhere in the critical strip off the critical line. 

We know that,
$$f(s)=-\left(\frac{1}{s}+J(s)\right)=\sum\limits_{m=1}^{\infty} T_{m}(s)=\sum\limits_{m=1}^{\infty}\frac{C_m}{\left(\frac{s}{2}\right)^{\overline{m}}}$$.

\begin{figure}[!htb]
\center   
\includegraphics[scale=.6]{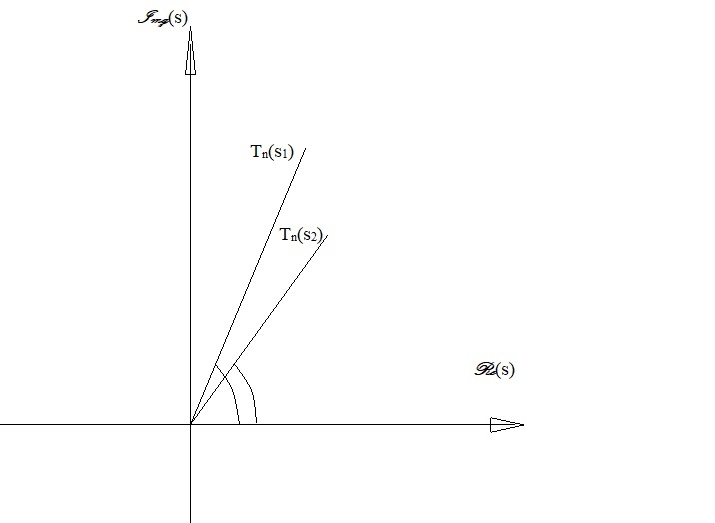}
\caption{Magnitude and argument (not necessarily less than $\pi/2$  as depicted) of the $nth$ series term $T_n(s_2)$ of $-f(s_2)$  is always less than the corresponding term $T_n(s_1)$ of $-f(s_1)$ if $\Re (s_1)< \Re (s_2)$ and $\Im (s_1) =\Im(s_2)$. }
\label{fig:image6}
\end{figure}

  Note that $-f(s)$ can be plotted by plotting $J(s)$ but  with the first term's coefficient $C_1$ taken as $ \left(\frac{1}{2}+\sum\limits_{n=1}^{\infty} e^{-n^2\pi}\right)$  instead of $\sum\limits_{n=1}^{\infty} e^{-n^2\pi}$.
  We, know already that from the center of the spiral the radial distance monotonously increases with each term such that it's maximum after counting all the terms i.e. from the spiral centre to the origin. The point Q is the center of the spiral. It has value $-f(s_1)$ which is reached  by laying the terms  $OA=T_1(s_1)$, $AB=T_2(s_1)$, $BC=T_3(s_1)$ and so on in the complex plane. We also draw concentric circles centered at the center of the spiral and having radius equal to the radial  distance from the center after the addition of each term (see fig. 6). Further, we decorate the corresponding terms $OA'=T_1(s_2)$, $AB'=T_2(s_2)$, $BC'=T_3(s_2)$ and so on of $-f(s_2)$ at points $ O, A, B, C....$ respectively. As a term with relatively smaller magnitude and argument replaces a term in the spiral,
  shorter radial distance covered by the spiral in that term. From the figure one can see that to draw the spiral geometry $-f(s_2)$ about point Q one has to move it's terms across the light red annular regions to move close to Q. So, the term OA' is moved to coincide point A' and point A across the last light red circular band. Similarly, second term AB' has to be moved across the second annular region (in light red) to coincide points B and B'. Repeating the same for rest of the terms to form a continuous spiral geometry $-f(s_{2})$. Thus, to draw the spiral geometry for $-f(s_2)$ about Q all the light red coloured annular regions are eliminated. This implies that the radial distance covered by each term  of $-f(s_2)$ (in yellow) about Q is less than that for $-f(S_1)$.
  Thus, we conclude $|\sum\limits_{i=n_1}^{n_2}T_n(s_1)|>|\sum\limits_{i=n_1}^{n_2}T_n(s_2)|$ and therefore,  $|f(s_{1})|>|f(s_{2})|$.
\begin{figure}[!htb]
\includegraphics[scale=.6]{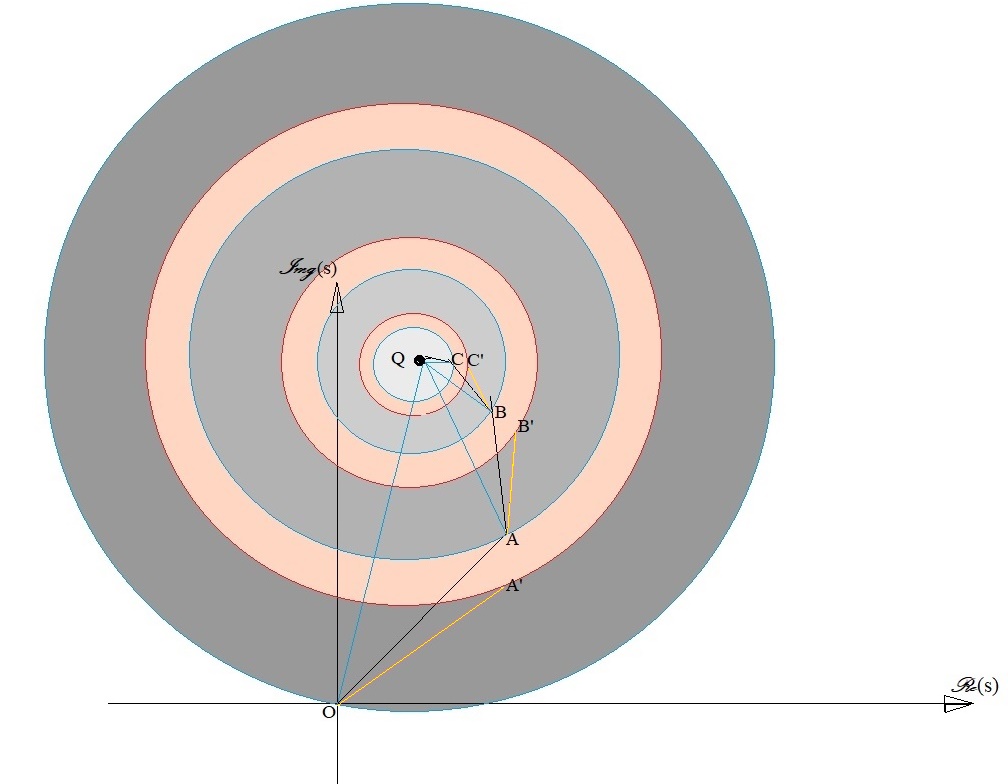}
\caption{Spiral geometry O-A-B-C-..-Q is term by term plot of $-f(s_{1})$ and the corresponding terms OA', AB', BC',... of $-f(s_{2})$ (in yellow) are placed from first term till the last term. The light red annular regions correspond (not equal to) the term wise difference in radial distance about Q, overall resulting in $|f(s_{1})|>|f(s_{2})|$.Given, $\Re (s_1)<\Re(s_2)$ and $\Im(s_1)=\Im(s_2)$.}
\label{fig:image6}
\end{figure}
If for two variables $s_{1}$ and $s_{2}$, where $s_1=\frac{1}{2}-\alpha-i|t|$ and $s_2=1-s_1^*=\frac{1}{2}+\alpha - i|t|$ where $\alpha \in (0,\frac{1}{2})$ and $t \in \mathbb{R}$  then as shown in the fig. 5,
 $$\bigg|T_{n}(s_{1})\bigg|>\bigg|T_{n}(s_{2})\bigg|$$.
$$\arg\left(T_{n}(s_{1})\right)>\arg\left(T_{n}(s_{2})\right)$$.

Thus, $$|f(\frac{1}{2}-\alpha - i|t|)|>|f(\frac{1}{2}+\alpha - i|t|)|=|f(\frac{1}{2}+\alpha + i|t|)|$$.
 
 This means nowhere in the critical strip can $f(s)$ be equal to $-f(1-s)$ off the critical line.
Here we can see that the non-trivial zeros correspond to purely imaginary values of $f(\frac{1}{2}+it)$ such that $f(s)+f(1-s)= f(s)+f(s^*)=0$.

Since, both the cases for falsification of the hypothesis are non-existent.
\textbf{Thus, the Riemann hypothesis is true.}\\

\section*{APPENDIX: On convergence of $J(s)$ }

In this section it's shown that the series function $J(s)$ converges in the critical strip region for $\vert \frac{s}{2}\vert\; \geq \; 1$. This is essentially the region of interest as in the region under the disk $\vert \frac{s}{2} \vert\;<\ 1$, it's known that no zero exists. \\
\;\;We use the Taylor series expansion of a new series $\Lambda(s)$. Using the integration test for convergence of $\Lambda(s)$ in the critical region and comparing the magnitude of coefficients of corresponding terms in $J(s)$, it's shown that $\vert J(s)\vert <\vert \Lambda(\frac{s}{2})\vert$. Thus, $\;J(s)\;$ converges.\\\\
Let  
$\Lambda (x)=\sum\limits_{n=1}^{\infty} e^{-n^{2}\pi x}$\\
Then, Taylor series expansion of $\Lambda(x)$ about point $x=1$ is\\

$\Lambda (x) = \sum\limits_{n=1}^{\infty} e^{-n^{2}\pi} +\left( 1-x\right)\sum\limits_{n=1}^{\infty} e^{-n^{2}\pi}\left(n^{2}\pi\right)+\frac{\left( 1-x\right)^2}{2}\sum\limits_{n=1}^{\infty} e^{-n^{2}\pi}\left(n^{2}\pi\right)^2 +... $\\

Also notice that all terms are positive for $0 < x < 1$.\\\\
Now, using the integration test for convergence  of $\Lambda(x)$ \cite{4},\\

$ a_{1}\geq\lim_{n\to\infty} \left( S_{n}\left(x\right)- I_{n}\left(x\right)\right)\geq 0 \;\;\;\;\;\;\;\;\;\;\;\;\;\;\;\;\;\;\;\;\;\;\;\;\;\;\;\;\;\;\;\;\;$ (1)\\
 here,
$ S_{n}\left(x\right)=\sum\limits_{j=1}^{n} e^{-j^{2}\pi x}$\\
and $a_{1} = e^{-\pi x}$\\

$ \lim_{n\to\infty} I_{n}\left(x\right)= \int\limits_{1}^{\infty} e^{-k^{2}\pi x} dk=\frac{1}{2\sqrt{x}}erfc(1)$\\

Where, erfc(x)  is the complimentary error function defined as,
$ erfc(x)=\frac{2}{\sqrt{\pi}}\int\limits_{x}^{\infty} e^{-t^{2}} dt$\\
 thus, from (1) we get\\

$ e^{-\pi x}+\frac{1}{2\sqrt{x}}erfc(1)\geq \Lambda(x) \geq \frac{1}{2\sqrt{x}}erfc(1) $\\\\
Thus, $\Lambda(s)$ is convergent only for $x>0$\\  

on expanding $\Lambda(x)$ from a real function to the complex plane\\

$\Lambda (s) = \sum\limits_{n=1}^{\infty} e^{-n^{2}\pi} +\left( 1-s\right)\sum\limits_{n=1}^{\infty} e^{-n^{2}\pi}\left(n^{2}\pi\right)+\frac{\left( 1-s\right)^2}{2}\sum\limits_{n=1}^{\infty} e^{-n^{2}\pi}\left(n^{2}\pi\right)^2 +... $\\
\\Where $s=\sigma +it\;\;\;\;\;$, $\sigma\;\in (0,1)$ and $ t \in \mathbb{R}$\\

but,\\

$ \vert \sum\limits_{n=1}^{\infty} e ^{-n^2\pi(\sigma + i t)}\vert \leq \vert\sum\limits_{n=1}^{\infty}e^{-n^{2} \pi(\sigma+ i t)}\vert =\sum\limits_{n=1}^{\infty} e^{-n^{2}\pi\sigma}\leq e^{-\pi\sigma} +\frac{1}{2\sqrt{\sigma}}erfc(1)$\\
thus,\\
$\Lambda (s)=\sum\limits_{n=1}^{\infty} e^{-n^{2}\pi s}$\\ is convergent for $\sigma >0$\\

Coming back to $J(s)$,\\

$J(s):= \sum\limits_{m=1}^{\infty}\frac{\sum\limits_{n=1}^{\infty} e^{-n^2\pi}\left(n^2\pi\right)^{m-1}}{\left(\frac{s}{2}\right)^{\overline{m}}}
$\\\\
Since, $\;\forall\; s\;$ in the critical strip :$ \vert s\vert\geq1$ and $k \in \mathbb{N}^+$, one can say that\\\\
$\;\;\;\;\;\vert s\vert\vert s + k\vert > k$\\\\
$\Rightarrow \frac{\vert s\vert}{k}>\frac{1}{\vert s + k\vert}$\\

$ \Rightarrow\frac{\vert\left(\frac{s}{2}\right)^k\vert}{k!}\;=\;\frac{\vert\frac{s}{2}\vert^{k}}{k!} \;>\;\frac{1}{\big| \left(\frac{s}{2}\right)\big|\big|\left(\frac{s}{2}+1\right)\big|...\big|\left(\frac{s}{2}+k-1\right)\big|\big|\left(\frac{s}{2}+k\right)\big|}\;=\; \frac{1}{\big|\left(\frac{s}{2}\right)\left(\frac{s}{2}+1\right)...\left(\frac{s}{2}+k-1\right)\left(\frac{s}{2}+k\right)\big|}  $\\

$\Rightarrow\vert \Lambda \left(\frac{s}{2}\right)\vert> \vert J(s)\vert\;\;\; \forall\; s$  in the critical strip outside the disk $\vert \frac{s}{2} \vert \leq 1 $ centered at origin. \\

Thus,

$J(s):= \sum\limits_{m=1}^{\infty}\frac{\sum\limits_{n=1}^{\infty} e^{-n^2\pi}\left(n^2\pi\right)^{m-1}}{\left(\frac{s}{2}\right)^{\overline{m}}}\;\;$
converges in our region of interest. \\\\Hence, proved!
\\\\
Therefore, the ratio of it's consecutive terms,\\
$ \frac{\vert T_{m+1}\vert}{\vert T_{m}\vert}= \frac{\sum\limits_{n=1}^{\infty} e^{-n^2\pi}\left(n^2\pi\right)^{m}}{\vert\frac{s}{2} + m\vert \sum\limits_{n=1}^{\infty} e^{-n^2\pi}\left(n^2\pi\right)^{m-1}}<1$.  

Another way to conclude the convergence of $J(s)$ in the critical strip is to use the fact that the Functional equation of zeta function converges in the critical strip and therefore all its component functions must.
	

\end{document}